\documentclass[10pt,twoside]{article}
\usepackage{amsfonts}
\usepackage{amssymb}
\usepackage{epsfig}
\usepackage{qtree}
\input xypic
\xyoption{all}
\LaTeXdiagrams
\xyoption{2cell}
\UseAllTwocells

\newtheorem{teo}{Theorem}

\newtheorem{cor}{Corollary}

\def\be{\begin{equation}}
\def\ee{\end{equation}}
\def\raa{\rightarrow}
\def\lra{\longrightarrow}

\def\A{{\mathbb A}}
\def\B{{\mathbb B}}
\def\C{{\mathbb C}}
\def\D{{\mathbb D}}

\def\circt{\ddot{\circ}}
\def\N{{\mathbb N}}
\def\Q{{\mathbb Q}}
\def\R{{\mathbb R}}
\def\Z{{\mathbb Z}}
\def\b{\fbox{\bf B}}
\def\c{\fbox{\bf C}}
\def\r{\fbox{\bf R}}
\def\m{\fbox{\bf M}}
\def\p{\fbox{\bf P}}
\def\boxt{\fbox{\bf B'}}
\def\boxcirct{\fbox{\bf C'}}

\def\PRdesc{{\bf PRdesc}}
\def\PRalg{{\bf PRalg}}
\def\PRfunc{{\bf PRfunc}}
\def\cat{{\bf Cat}}
\def\catxn{{\bf CatXN}}

\def\la{\langle}
\def\ra{\rangle}
\newenvironment{examp}{ \stepcounter{examnum} {\bf \noindent Example:}}{$\Box$}

\newcounter{examnum}[section]
\newcounter{remarnum}[section]
\begin{document}
\title{Towards a Definition of an Algorithm}
\author{Noson S. Yanofsky\footnote{Department of Computer and Information Science,
Brooklyn College, CUNY,
Brooklyn, N.Y. 11210. And Computer Science Department,
The Graduate Center, CUNY,
New York, N.Y. 10016.
e-mail: noson@sci.brooklyn.cuny.edu}}
\maketitle
\begin{abstract}
\noindent We define an algorithm to be the set of programs that
implement or express that algorithm. The set of all programs is
partitioned into equivalence classes. Two programs are equivalent if
they are essentially the same program. The set of equivalence
classes forms the category of algorithms. Although the set of
programs does not even form a category, the set of algorithms form a
category with extra structure. The conditions we give that describe
when two programs are essentially the same turn out to be coherence
relations that enrich the category of algorithms with extra
structure. Universal properties of the category of algorithms are
proved.

\vspace{.3 in}
\end{abstract}

\section{Introduction}
In their excellent text {\it Introduction to Algorithms, Second Edition }\cite{Corman}, Corman,
Leiserson, Rivest, and Stein begin Section 1.1 with a definition of an algorithm:
\begin{quote}
Informally, an {\bf algorithm} is any well-defined computational procedure that takes some value, or set
of values, as {\bf input} and produces some value, or set of values, as {\bf output}.
\end{quote}
Three questions spring forward:
\begin{enumerate}
\item ``Informally''? Can such a comprehensive and highly technical book of 1180 pages not have a ``formal''
definition of an algorithm?
\item What is meant by ``well-defined?''
\item The term ``procedure'' is as vague as the term ``algorithm.'' What is a ``procedure?''
\end{enumerate}
Knuth \cite{Knuth,Knuth2} has been a little more precise in specifying the requirements demanded for an algorithm.
But he writes ``Of course if I am pinned down and asked to explain more precisely what I mean by these remarks,
I am forced to admit that I don't know any way to define any particular algorithm except in a programming language.''
(\cite{Knuth2}, page 1.)

Although algorithms are hard to define, they are nevertheless real mathematical objects. We name and talk
about algorithms with
phrases like ``Mergesort runs in $n$ lg $n$ time''. We quantify over all algorithms, e.g.,
``There does not exist an algorithm to solve the halting problem.'' They are as ``real'' as the number $e$
or the set $\Z$. See \cite{Dean} for an excellent philosophical overview of the subject.

Many researchers have given definitions over the years. (Refer to \cite{BlassGur} for a historical
survey of some of these definitions. One must also read the important works of Yiannis Moschovakis, e.g., \cite{Mosch}.)
Many of the given definitions are of the form ``An algorithm is a program in this language/system/machine.''
This does not really conform to the current usage of the word ``algorithm.'' Rather, this is more in tune
with the modern usage of the word ``program.'' They all have a feel of being a specific implementation
of an algorithm on a specific system. Imagine a professor teaching a certain algorithm to a class and then
assigning the class to go home and program the algorithm. In any class with the moral abhorrence of cheating,
the students will return many
{\it different} programs implementing the {\it same} algorithm. We would not call each of these different programs an algorithm.
Rather the different programs are implementations of a single algorithm. And yet some researcher do call each of those
programs a different algorithm, e.g. \cite{BlassDerGur}.
We would like to propose another definition.

Consider Figure 1.

\begin{figure}[!h]
\centering
  \includegraphics[height=3.25 in]{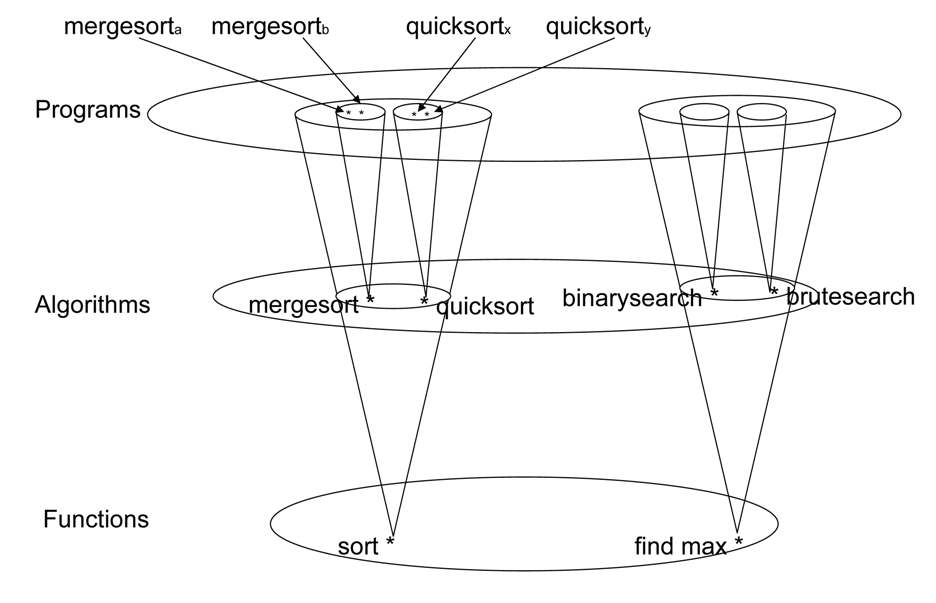}\\
\caption{Programs, Algorithms and Functions.}
\end{figure}

At the bottom of the figure is the set of all functions. Two
functions are highlighted: the sort function and the function that
outputs the maximum of its inputs. On top of the figure is the set
of all programs. For every function there is a set of programs that
implement that function. We have highlighted four programs that
implement the sort function: $\mathrm{\bf mergesort_a}$ and
$\mathrm{\bf mergesort_b}$ are two different programs that implement
the algorithm mergesort. Similarly $\mathrm{\bf quicksort_x}$ and
$\mathrm{\bf quicksort_y}$ are two different implementations of the
algorithm quicksort. There are also many different programs that
implement the max function. $\mathrm{\bf mergesort_a}$ and
$\mathrm{\bf mergesort_b}$ are grouped in one subset of all the
programs that implement the sort function. This subset will
correspond to the mergesort algorithm. Similarly, $\mathrm{\bf
quicksort_x}$ and $\mathrm{\bf quicksort_y}$ are grouped together
and will correspond to the quicksort algorithm. There are similar
groupings for a binary search algorithm that finds the max of a list
of elements. There are also other algorithms that find the max. This
intuition propels us to define an algorithm as the set of all
programs that implement the algorithm.

We define an algorithm analogously to the way that Gottlob Frege
defined a natural number. Basically Frege  says that the number 42
is the equivalence class of all sets of size 42. He looks at the
conglomerate of all finite sets and makes an equivalence relation.
Two finite sets are equivalent if there is a one-to-one onto
function from one set to the other. The set of all equivalence
classes under this equivalence relation forms the set of natural
numbers. For us, an algorithm is an equivalence class of programs.
Two programs are part of the same equivalence class if they are
``essentially'' the same. Each program is an expression (or an
implementation) of the algorithm, just as every set of size 42 is an
expression of the number 42.

For us, an algorithm is the sum total of all the programs that express it. In other
words, we look at all computer programs and partition them into different subsets. Two programs
in the same subset will be two implementations of the same algorithm. These two programs are ``essentially''
the same.

What does it mean for two programs to be ``essentially'' the same? Some examples are in order:
\begin{itemize}
\item One program might perform $Process_1$ first and then perform an unrelated $Process_2$ after.
The other program will perform the two unrelated processes in the opposite order.
\item One program might perform a certain process in a
loop $n$ times and the other program will unwind the loop and perform it $n-1$ times and then perform the the process again
outside the loop.
\item One program might perform two unrelated processes in one loop, and the other
program might perform each of these two processes in its own loops.
\end{itemize}
In all these examples, the two programs are
definitely performing the same function, and everyone would agree that both programs are implementations of the same
algorithm. We are taking that subset of programs to be the definition of an algorithm.

Many relations that say when two programs are essentially the same will be given. However, it is doubtful
that we have the final word on this. Hence the word ``Towards'' in the title.
Whether or not two programs are essentially the same, or whether or not a program is
an implementation of a particular algorithm is really a subjective decision. Different relations can be given for different purposes.
We give relations that most people can agree on that these two programs are essentially the same, but we are well aware of the
fact that others can come along and give more, less or different relations. The important realization is that the relations that we feel are the most
obvious turn out to be relations that correspond to standard categorical coherence rules. When we mod-out by any set of relations, we get more structure.
When we mod-out by these
relations, our set of programs become a category with more structure. Our goal is not to give the final word on the topic, but
to point out that this is a valid definition of an algorithm and that the equivalence classes of algorithms has more
structure than the set of programs.

We consider the set of all programs which we might call $\mbox{\bf Programs}$. An equivalence relation $\approx$ of
``essentially the sameness''
is then defined on this set. The set of equivalence classes $\mbox{\bf Programs}/\approx$ shall
then be called $\mbox{\bf Algorithms}$. There is a nice onto function from $\phi:\mbox{\bf Programs} \lra
\mbox{\bf Algorithms}$, that takes every program $P$ to the equivalence class $\phi(P)= \left[ P \right]$.
One might think of
any function $\psi:\mbox{\bf Algorithms} \raa\mbox{\bf Programs}$ such that $\phi \circ \psi=Id_{\mbox{\bf Algorithms}}$
as an ``implementer.'' $\psi$ takes an algorithm to an implementation of that algorithm.

To continue with this line of reasoning, there are many different algorithms that perform the same function.
For example, Kruskal's algorithm and Prim's algorithm are two different ways of finding a minimum spanning tree
of a weighted graph. Quicksort and Mergesort are two different algorithms to sort a list. There exists
an equivalence relation on the set of all algorithms. Two algorithms are equivalent $\approx'$ if they
perform the same function. We obtain $\mbox{\bf  Algorithms} /\approx'$ which we might call $\mbox{\bf Comp. Functions}$
or computable functions.
It is an undecidable problem to determine when two programs perform the same computable function.
Hence we might not be able to effectively
give the relation $\approx'$, nevertheless it exists. Even if we were able to give the relation, that would not mean that
the word problem (i.e., telling when two different equivalence classes of descriptions are equivalent) is solvable.
Nevertheless, there is an onto function $\phi':\mbox{\bf Algorithms} \lra
\mbox{\bf Comp. Functions}$.

We summarize our intentions with the following picture.
$$\xymatrix{
\mbox{Programing} &\mbox{Computer Science}&\mbox{Mathematics}\\
\mbox{\bf Programs} \ar@{>>}[r]& \mbox{\bf Algorithms} \ar@{>>}[r]&
\mbox{\bf Comp. Functions }\\
}$$
Programs are what programmers, or
software engineers deal with. Algorithms are the domain of computer scientists.
Computable functions are of interest to pure mathematicians.

With this picture in mind, we can explain  other equivalence relations describing program
``sameness''. One can give many different equivalence relations but they must fall within the two extremes.
One extreme says that no two programs are really the same, i.e., every program is essentially an
algorithm. In that case {\bf Programs} = {\bf Algorithms}. This extreme case is taken up by \cite{BlassDerGur}.
In contrast, another extreme is to say
that two programs are the same if they perform the same operation or are bisimilar.
In that case {\bf Algorithms} = {\bf Comp. Functions}. In this paper we choose a middle way.
Others can have other equivalence relations but they must fall in the middle. There are finer and courser
equivalence relations than ours. There will also be unrelated equivalence relations. For every
equivalence relation, the set of algorithms will have a particular structure.

In our scheme, {\bf Programs} will form a directed graph with a composition of arrows and a distinguished loop on
every vertex. However they will not have the structure of
a true category: the composition will not be associative and the distinguished loops will not act like the identity.
In contrast, {\bf Algorithms} will be a real category with extra structure: a Cartesian product structure and a weak parameterized natural
number object (a categorical way of saying that the category
is closed under recursion). This category will turn out to be an initial category in the 2-category of
all categories with products and weak parameterized natural number objects.

Others have studied similar categories before.  Joyal in an
unpublished manuscript about ``arithmetical universes''. (see
\cite{Maietti} for a history) as well as \cite{Burroni},
\cite{Roman} and \cite{Pareandroman} have looked at the free
category with products and a {\it strong} natural number object.
Marie-France Thibault \cite{mft} has looked at a Cartesian {\it
closed} category with a weak natural number object. They
characterized what type of functions can be represented in such
categories. Although related categories have been studied, the
connection with the notion of an algorithm has never been seen. Nor
has this category ever been constructed as a quotient of a
syntactical graph.

We are not trying to make any ontological statement about the existence of algorithms. We are merely
giving a mathematical way of describing how one might think of an algorithm. Human beings dealt with
rational numbers for millennia before mathematicians decided that rational numbers are equivalence
classes of pairs of integers:
$$\Q\quad =\quad \{(m,n)\in \Z \times \Z|n\neq 0 \} / \approx $$
where
$$(m,n)\approx (m',n') \mbox{ iff } mn'=nm'.$$
Similarly, one can think of the existence of algorithms in any way that one chooses. We are simply
offering a mathematical way of presenting them.

There is a interesting analogy between thinking of a rational number as an equivalence class of pairs of integers
and our definition of an algorithm as an equivalence class of programs. Just as a rational number can only
be expressed by an element of the equivalence class, so too, an algorithm can only be expressed by presenting
an element of the equivalence class. When we write an algorithm, we are really writing a program. This explains the
quote from
Knuth's given in the beginning of this paper. Pseudo-code
is used to allow for ambiguities and not show any preference for a language. But it is, nevertheless, a program.

Another applicable analogy is just as a rational number by
itself has no structure (it is simply an equivalence class of pairs
of integers), so too, an algorithm has no structure. In contrast, the
set of rational numbers has much structure. So too, the set
(category) of algorithms has much structure. $\Q$ is the smallest
field that contains the natural numbers. We shall see in Section 4
that the category of algorithms is an initial category with a
product and a weak natural number object.

When a human being talks about a rational number, he prefers to use the pair $(3,5)=\frac{3}{5}$ as opposed to
the equivalent pair $(6,10)$, or
the equivalent $(3000, 5000)$. One might say that the rational number $(3,5)$ is a ``canonical representation'' of the
equivalence class to which it belongs. It would be nice if there was a ``canonical representation'' of an
algorithm. We speculate further on this ideas in the last section of this paper.

The question arises as to which programming language should we use?
Rather than choosing one programming language to the exclusion of others, we look at a
language of descriptions of primitive recursive functions. We choose this language because
of its beauty, its simplicity of presentation, and the fact that most readers can easily become familiar with this language.
The language of descriptions of primitive recursive functions basically has three operations: Composition, Bracket, and Recursion.
A primitive recursive function
can be described in many different ways. A description of a primitive recursive function is basically the same
thing as a program in that it tells how to calculate a function. There is a basic correlation
between programming concepts and
the operations in generating descriptions of primitive recursive functions: recursion is like a loop, composition is
sequential processing, and bracket is parallel processing. We are well aware that we are limiting ourselves because the set
of primitive recursive functions is a proper subset of the set of all computable functions. By limiting
ourselves, we are going to get a proper subset of all algorithms. Even though we are, for the present time, restricting
ourselves, we feel that the results we will get are interesting in their own right. There is an ongoing project to extend this
work to all recursive functions \cite{ManinYano}.

\vspace{.5in}

There is another way to view this entire endeavor. What we are creating here is an operad.
Operads are a universal algebraic/categorical way of describing extra algebraic structure.
Recently operads have become very popular with algebraic topologists and
people who study quantum field theories. We are creating an operad that describes some of the extra
structure that exists on the set of total functions of a certain type. With such total functions
one can compose, do recursion, and take the product of those functions. We than can look
at the algebra of this operad generated by all total functions from powers of $\N$ to $\N$.
One then can examine the subalgebra generated by basic or initial functions (this essentially
is our $\PRdesc$). We go further and look at a quotient of this subalgebra by using more relations
(this essentially is our $\PRalg$). We show in section 4 of this paper that this quotient subalgebra is
an initial object in a certain 2-category.
This operadic viewpoint is further elaborated and used in \cite{ManinYano} where we tackle the harder problem
of all recursive functions.

\vspace{.5in}

There is a fascinating correspondence between this work and similar work in low-dimensional topology and related work
in topological quantum field theory (TQFT). This correspondence is in the spirit of \cite{BaezStay} and \cite{Abramsky}
where they show that using the powerful language of category theory there are many similar phenomena in
low-dimensional topology, quantum physics, and logic. In order for us to express this correspondence, we are going to
have to assume  some knowledge of the basic yoga of low-dimensional topology. If this is not known, then simply
skip this paragraph. For clarity's sake, we shall concentrate on the category of braids. However, we could have described similar
correspondences with tangles, ribbons, cobordisms, etc. Similar to our three levels of structure,
$$\xymatrix{
\mbox{\bf Programs} \ar@{>>}[r]& \mbox{\bf Algorithms} \ar@{>>}[r]&
\mbox{\bf Comp. Functions }\\
}$$
there are three levels of objects in low-dimensional topology:
$$\xymatrix{
\mbox{\bf Braid Projections} \ar@{>>}[r]& \mbox{\bf Braid Groups} \ar@{>>}[r]&
\mbox{\bf Symmetric Groups. }\\
}$$
With these, there are the following analogies.
\begin{itemize}
\item Just as we can only represent an algorithm by giving a program, so too, the only way to represent a braid is by giving
a braid projection.
\item Just as our set of {\bf Programs} does not have enough structure to form a category, so too, the set of
{\bf Braid Projections} does not have a worthwhile structure. One can compose braid projections sequentially and parallel. But there
is no associativity. There are identity braids, but when sequentially composed with other
braid projections, they do not act like projections. There are inverse braid projections, but when sequentially composed with
the original projection, there is no identity projection.
\item Just as we can get the category of algorithms by looking at equivalence classes of programs, so too, we can get braids
by looking at equivalence classes of braid projections. With braid projections we look at Reidermeister moves to determine when
two braid projections are really the same. Here we look at relations stated in this paper to tell when two programs are the same.
\item Just as we are not giving the final word about what relations to use, so too, there is no final word about which
Reidermeister moves to use. Depending on your choice, you will get
braids, ribbons, oriented ribbons etc.
\item Just as our category of {\bf Algorithms} is the free category with
products and a weak natural number object generated by the empty
category, so too, the category of {\bf Braids} is the free braided
monoidal category generated by one object.
\item Just as we can go down to the level of functions by making two algorithms that perform the same function equivalent, so to, we can
add a relation that two strings can cross each other and get the {\bf Symmetric Groups}.
\item Just as the main focus of computer scientists are algorithms and not programs, so to, the main focus of topologists is braids and
not braid diagrams.

\end{itemize}
There is obviously much more to explore with these analogies. There also should be a closer relationship between these fields. After all,
some of our relations are very similar to Reidermeister moves.

\vspace{.5in}

Section 2 will review the basics of primitive recursive functions and show how they may be described
by special labeled binary trees. Section 3 will then give many of the relations that tell when two descriptions
of a primitive recursive function are ``essentially'' the same. Section 4 will discuss the structure of the set of all algorithms. We shall
give a universal categorical description of the category of algorithms. This is the only Section that uses
category theory in a non-trivial way. Section 5 will
discuss complexity results and show how complexity theory fits into our framework. We conclude this
paper with a list of possible ways this work can progress.

At this point it is appropriate to say what this paper {\it is not}.
\begin{itemize}
\item We have no ambition to say anything new about primitive recursive functions. We are only using
descriptions of primitive recursive functions as a simple programming language with three operations. Nor are we
saying anything about a relationship between programming languages and primitive recursive functions.
\item Nothing new will be said about category theory. Rather, we are making a link of these categories and the concept of an algorithm.
\item We will not say anything new about program semantics. Our equivalence relations are between descriptions that
correspond to the same function.
\end{itemize}

Rather, what we are doing here is giving a novel definition of an algorithm and showing that the
the set of all algorithms has more manageable structure than the set of all programs. We are also showing
that categorical coherence relations correspond to rules saying when two programs are essentially the same.

\vspace{.5in}

Yuri Manin has incorporated an earlier draft \cite{Yano}
of this paper into his second edition of his {\it A Course in Mathematical Logic} \cite{Manin}. Within Chapter
IX of that book he describes the constructions given in this paper using the language of PROPs and operads that
are of interest to mathematicians and theoretical physicists. This earlier draft \cite{Yano} was also discussed in \cite{BlassDerGur}.

\vspace{.5in}

\noindent {\bf Acknowledgement.} Alex Heller (OBM), Florian Lengyel,
and Dustin Mulcahey were forced to listen to me working through this
material. I am indebted to them for being total gentlemen while
suffering in silence. I am grateful to Rohit Parikh, Karen Kletter,
Walter Dean, and the entire Graduate Center Logic and Games Seminar
gang for listening to and commenting on this material. Shira
Abraham, Ximo Diaz Boils, Juan B. Climent, Eva Cogan, Joel Kammet,
and Matthew K. Meyer read through and commented on many versions of
this work. This paper would have been, no doubt, much better had I
listened to all their advice. I owe much to Yuri Manin and Grigori
Mints for taking an interest in this work and for long e-mail
discussions. This work was inspired by a talk that Yuri Gurevich
gave on his definition of an algorithm.

\section{Descriptions of Primitive Recursive Functions}
Rather than talking of computer programs, {\it per se}, we shall talk of descriptions of
primitive recursive functions. For every primitive recursive function, there are many different
methods of ``building-up'', or constructing the function from the basic functions. Each method is similar to a program.

We remind the reader that primitive recursive functions $\N^n \raa\N$ are ``basic'' or ``initial'' functions:
\begin{itemize}
\item null function $n:\N \raa\N$ where $n(x)=0$
\item successor function $s:\N \raa\N$ where $s(x)=x+1$
\item for each $k \in \N$ and for each $i\leq k$, a projection function $\pi^k_i:\N^k \raa\N$ where
$\pi^k_i(x_1,x_2, \ldots x_k) = x_i$
\end{itemize} and functions constructed from basic functions through a finite number of compositions and recursions.

We shall extend this definition in two non-essential ways. An $n-$tuple  of primitive recursive functions
$(f_1,f_2, \ldots f_n):\N^m \raa\N^n$,
shall also be called a primitive recursive function. Also, a constant function
$k: \ast \raa\N$ is called a primitive recursive function because for every $k\in \N$, the constant
map may be written as $s\circ s \circ \cdots \circ s \circ n$.

Let us spend a few minutes reminding ourselves of basic facts about recursion. The simplest form
of recursion is
for a given integer $k$ and a function $g:\N \raa\N$. From this one constructs $h:\N \raa\N$
as follows

 $h(0)=k$

 $h(n+1)=g(h(n)).$

A more complicated form of recursion --- and the one we shall employ ---
 is for a given function $f:\N^k \raa\N^m$ and a given function
$g:\N^k \times \N^m  \raa\N^m$. From this one constructs $h:\N^k \times \N \raa\N^m$ as

 $h(x,0)=f(x)$

 $h(x,n+1)=g(x,h(x,n))$

\noindent where $x \in \N^k$ and $n \in \N$.

The most general form of recursion, and the definition usually given for
primitive recursive functions is for a given function $f:\N^k \raa\N^m$ and a given function
$g:\N^k \times \N^m \times \N \raa\N^m$. From this, one constructs $h:\N^k \times \N \raa\N^m$

 $h(x,0)=f(x)$

 $h(x,n+1)=g(x,h(x,n),n)$

\noindent where $x \in \N^k$ and $n \in \N$.

We shall use the middle definition of recursion because the extra input variable in $g$ does not add anything
\cite{Gladstone}. It simply  makes things
unnecessarily complicated. However, we are certain that any proposition that can be said about
the second type of recursion, can also be said for the third type. See \cite{TTT1} Section 7.5, and \cite{TTT} Section 5.5.

Although primitive recursive functions are usually described as closed only under composition and recursion,
there is, in fact, another implicit operation for which the functions are closed: bracket. Given
primitive recursive functions $f:\N^k \raa\N$ and $g:\N^k \raa\N$, there is a primitive recursive function
$h=\la f,g \ra :\N^k \raa\N \times \N$. $h$ is defined as $$h(x)= (f(x),g(x))$$ for any $x \in \N^k$. We shall
see that having this bracket operation is almost the same as having a product operation.

In order to save the eyesight of our poor reader, rather than writing too
many exponents, we shall write a power of the set $\N$
for some fixed but arbitrary number as $\A, \B, \C$ etc.
With this notation, we may write the recursion operation as follows: from functions
$f:\A \raa\B$ and $g:\A \times \B \raa\B$ one constructs $h:\A \times \N \raa\B$.

If $f$ and $g$ are functions with the appropriate source and targets, then we shall write
their composition as $h=f\ \circ g$. If they have the appropriate source and target for the
bracket operations, we shall write the bracket operation as $h=\la f,g \ra$. We are in need of a similar
notation for recursion. So if there are $f:\A \raa\B$ and $g:\A \times \B \raa\B$ we shall write the
function that one obtains from them through recursion as $h=f \sharp g :\A \times \N \raa\B$

We are going to form a directed graph that contains all the descriptions of primitive recursive functions. We shall call
this graph $\PRdesc$. The vertices of the graph shall be powers of the natural number
$\N^0=\ast, \N, \N^2, \N^3, \ldots $. The edges of the graph shall be descriptions of
primitive recursive functions. One should keep in mind the following intuitive picture.

\vspace{.2in}

$$\xymatrix{ \cdots \ar[r] \ar@/^/[r] \ar@/_/[r]& \N^k \ar@/^.3in/[l]\ar@/^/[rrr] \ar@/_/[rr]
\ar[r] \ar@/^/[r] \ar@/_/[r] & \cdots \ar@/_.1in/[l]\ar[r] \ar@/^/[r] \ar@/_/[r]
&  \N^4 \ar[r] \ar@/^/[r] \ar@/_/[r]&\N^3 \ar[r] \ar@/^/[r] \ar@/_/[r]& \N^2 \ar[r] \ar@/^/[r] \ar@/_/[r]\ar@/_.1in/[lll]
\ar@/^.3in/[l]& \N \ar@/_.3in/[l]\ar@/^.1in/[l]&
\ast \ar[l] \ar@/^/[l] \ar@/_/[l].}$$

\vspace{.2in}

\subsection{Trees}
Each edge in $\PRdesc$ shall be a labeled
binary tree whose leaves are basic functions and whose internal nodes are labeled by {\bf C, R} or {\bf B}
for composition, recursion and bracket. Every internal node of the tree shall be derived from its
left child and its right child. We shall use the following notation:
\qtreecenterfalse
\Tree[{$f:\A\raa \B$} {$g:\B\raa \C$} ] .{$g\circ f:\A \raa \C$\\
\c} \qquad  \Tree[ {$f:\A\raa \B$} {$g:\A \times \B\raa \B$} ] .{$h=f\sharp
g:\A\times \N \raa \B$\\ \r} \qquad  \Tree[ {$f:\A\raa \B$} {$g:\A\raa
\C$} ] .{$\la f, g \ra :\A \raa \B \times \C$\\ \b}

$\PRdesc$ has more structure than a simple graph. There is a composition of edges. Given a tree $f:\A \raa\B$
and a tree $g:\B \raa\C$, there is another tree $g \circ f:\A \raa\C$. It is, however, important to
realize that $\PRdesc$ is {\it not} a category. For three composable edges, the trees
$h\circ ( g \circ f)$ and $(h \circ g) \circ f$ exist and they perform the same operation, but they are, nevertheless,
different programs and different trees. There is a composition of morphisms, but this composition is not
associative.

Furthermore, for each object $\A$ of the
graph, there is a distinguished morphism $\pi^\A_\A:\A \raa\A$ which does
not act like an identity. It is simply a function whose output is the same as its input.
\subsection{Some Macros}
Because the trees that we are going to construct can quickly become large and cumbersome, we will
employ several programming shortcuts, called {\it macros}. We use the macros to improve readability.

\noindent{\bf Multiple Projections.} There is a need to generalize the notion of a projection. The $\pi^k_i$
accept $k$ inputs and outputs one number. A multiple projection takes $k$ inputs and outputs $m$ outputs.
Consider $\A=\N^k$ and the sequence
$X= \la x_1,x_2, \ldots , x_m \ra$ where each $x_i$ is in $\{1,2, \ldots, k \}$. Let $\B=\N^m$, then for every $X$ there exists
$\pi^{\N^k}_{\N^m}=\pi^\A_\B:\A \raa\B$ as
$$ \pi^\A_\B = \la \pi^\A_{x_1}, \la \pi^\A_{x_2}, \la \ldots, \la \pi^\A_{x_{m-1}},\pi^\A_{x_m} \ra \ra \ldots \ra.$$
In other words, $\pi^\A_\B$ outputs the proper numbers in the order described by $X$. Whenever possible, we shall be
ambiguous with superscripts and subscripts.

Setting $$X=I=\la 1,2,3, \ldots, n \ra$$ we have what looks like the identity functions. Setting
$$X=\triangle=\la 1,2,3, \ldots, n, 1,2,3, \ldots, n \ra$$ we get the diagonal function.
\vspace{.2in}

\noindent {\bf Products.} We would like a product of two maps. Given $f:\A\raa\B$ and $g:\C\raa\D$, we would like
$f \times g:\A\times \C  \raa\B \times \D$. The product can be defined using the bracket as
$$f \times g = \la f \circ \pi^{\A \times \C}_\A, g \circ \pi^{\A \times \C}_\C \ra$$
or in terms of trees

\qquad \qquad \qquad \Tree [ {$f:\A\raa\B$} {$g:\C\raa\D$} ] .{$f \times g:\A \times \C  \raa\B \times \D$\\ \p}

\vspace{.2in}

is defined ($=$) as the tree

\vspace{.2in}

\Tree [ [ {$\pi^{\A \times \C}_\A:{\A \times \C}\raa\A $}  {$f:\A\raa\B$}
].{$f \circ \pi^{\A \times \C}_\A :\A\times \C  \raa\B$ \\ \c}
[ {$\pi^{\A \times \C}_\C:{\A \times \C}\raa\C $} {$g:\C\raa\D$}
].{$g \circ \pi^{\A \times \C}_\C :\A\times \C  \raa\D$ \\ \c}
].{$f \times g = \la f \circ \pi^{\A \times \C}_\A, g \circ \pi^{\A \times \C}_
\C \ra :\A\times \C  \raa\B \times \D$ \\ \b}
\vspace{.2in}

\noindent{\bf Diagonal Map.} A diagonal map will be used. A diagonal map is a map $\triangle:\A \raa\A \times \A$
where $x \mapsto (x,x)$. It can be defined as

$\triangle:\A \raa\A \times \A$ \qquad =
\Tree [{$\pi^\A_\A:\A \raa\A$} {$\pi^\A_\A:\A \raa\A$} ].{$\la \pi^\A_\A, \pi^\A_\A \ra: \A \raa\A \times \A$. \\ \b}

We took the bracket operation as fundamental and from the bracket operation we derived the product
operation and the diagonal map. We could have just as easily taken  the product and the diagonal as fundamental
and constructed the bracket as
$$\xymatrix{
\A \ar[rr]^{\la f,g \ra} \ar[rdd]_\triangle  && \B \times \C
\\
\\
& \A \times \A. \ar[ruu]_{f\times g}
}$$

\vspace{.2in}

\noindent{\bf Twist Map.} We shall need to switch the order of inputs and outputs. The twist map shall
be defined as
$$tw_{\A,\B} = \pi^{\A \times \B}_\B \times \pi^{\A \times \B}_\A :\A \times \B \raa\B \times \A.$$

Or in terms of trees:

$tw_{\A,\B}:\A \times \B \raa\B \times \A$
\qquad =
\Tree [{$\pi^{\A \times \B}_\B :\A \times \B \raa\B$} {$\pi^{\A \times \B}_\A :\A \times \B \raa\A$}
].{$\pi^{\A \times \B}_\B \times \pi^{\A \times \B}_\A :\A \times \B \raa\B \times \A$ \\ \p}
\vspace{.2in}

\noindent{\bf Second Variable Product.}
Given a function $g_1:\A \times \B \raa\B$ and a function $g_2:\A \times \B \raa\B$, we would like to take
the product of these two functions while keeping the first variable fixed. We define the operation
$$g_1 \boxtimes g_2:\A \times \B \times \B \raa\B \times \B$$
on elements as follows
$$(g_1 \boxtimes g_2)(a,b_1,b_2) = (g_1(a,b_1), g_2(a, b_2)).$$
In terms of maps, $\boxtimes$ may be defined from the composition of the following maps:
$$g_1 \boxtimes g_2= (g_1 \times g_2)\circ (\pi_\A \times tw_{\A,\B} \times \pi_\B) \circ (\triangle \times
\pi^{\B \times \B}_{\B \times \B}):$$
$$\A \times \B \times \B \raa\A \times \A \times \B \times \B \raa\A \times \B \times \A \times \B \raa\B \times \B.$$
Since the second variable product is related to the product which is derived from the bracket, we write it as

\qquad \qquad \qquad \qquad \Tree [{$g_1:\A \times \B \raa\B$} {$g_2:\A \times \B \raa\B$}
].{$g_1 \boxtimes g_2:\A \times \B \times \B \raa\B \times \B$\\ \boxt}
\vspace{.2in}

\noindent{\bf Second Variable Composition.}
Given a function $g_1:\A \times \D \raa\B$ and a function $g_2:\A \times \C \raa\D$, we would like to compose
the output of $g_2$ into the second variable of $g_1$. We define the operation
$$g_1 \circt g_2:\A \times \C \raa\B $$
on elements as follows
$$(g_1 \circt g_2)(a,c)= g_1(a,g_2(a,c)).$$
In terms of maps, $\circt$ may be defined as the composition of the following maps
$$g_1 \circt g_2= (g_1) \circ (\pi^\A_\A \times g_2) \circ (\triangle \times \pi^\C_\C): $$
$$\A \times \C \raa\A \times \A \times \C \raa\A \times \D \raa\B$$
We write second variable composition  as

\qquad \qquad \qquad \qquad\Tree [ {$g_2:\A \times \C \raa\D$} {$g_1:\A \times \D \raa\B$}
].{$g_1 \circt g_2:\A \times \C \raa\B $\\ \boxcirct}

\section{Relations}
Given the operations of composition, recursion and bracket, what does it mean for us to say that two
descriptions of a primitive recursive function are ``essentially'' the same? We shall examine
these operations, and give relations to describe when two trees are essentially the same. If two trees are
exactly alike except for a subtree that is equivalent to another tree, then we may replace the subtree
with the equivalent tree.

\subsection{Composition}

\noindent{\bf Composition is Associative.} That is, for any three composable maps $f$, $g$ and $h$, we have
$$h \circ (g \circ f) \approx (h \circ g) \circ f.$$
In terms of trees, we say that the following two trees are equivalent:

\qtreecenterfalse

\Tree [ [ {$f:\A\raa\B$} {$g:\B\raa\C$} ] .{$g\circ f:\A\raa\C$ \\ \c} {$h:\C \raa\D$}
].{$h\circ (g\circ f):\A\raa\D$\\ \c}
$\approx$
\Tree [ {$f:\A\raa\B$} [ {$g:\B\raa\C$} {$h:\C \raa\D$} ].{$h\circ g:\B\raa\D$\\ \c}
].{$(h\circ g) \circ f:\A\raa\D$\\ \c}

\vspace{.2in}

\noindent{\bf Projections as Identity of Composition.}
The projections $\pi^\A_\A$ and $\pi^\B_\B$ act like identity maps. That means for any $f:\A \raa\B$,
we have
$$f \circ \pi^\A_\A \approx f \approx \pi^\B_\B \circ f.$$
In terms of trees this amounts to

\Tree [{$\pi^\A_\A:\A \raa\A$} {$f:\A\raa\B$} ].{$f \circ \pi^\A_\A :\A \raa\B$\\ \c}
$\approx$
 {$\quad f:\A\raa\B \quad$}
$\approx$
\Tree[ {$f:\A\raa\B$} {$\pi^\B_\B:\B \raa\B$} ].{$\pi^\B_\B \circ f:\A \raa\B $\\ \c}

\vspace{.2in}

\noindent{\bf Composition and the Null Function.}
The null function always outputs a $0$ no matter what the input is. So for any function $f:\A\raa\N$, if
we are going to compose $f$ with the null function, then $f$ might as well be substituted with a projection, i.e.,
$$n \circ f \approx n \circ \pi^\A_\N.$$
In terms of trees:

\Tree [ [$f_1$ $f_2$ $\cdots$ $f_k$ ].{$f:\A\raa\N$} {$n:\N\raa\N$} ].{$n \circ f:\A\raa\N$\\ \c}
$\approx$
\Tree [{$\pi^\A_\N:\A\raa\N$} {$n:\N\raa\N$} ].{$n \circ \pi^\A_\N:\A \raa\N.$ \\ \c }

Notice that the left side of the left tree is essentially ``pruned.'' Although there is much information
on the left side of the left tree, it is not important. It can be substituted with another tree that
does not have that information.

\subsection{Composition and Bracket}
\noindent{\bf Composition Distributes Over the Bracket on the Right.}
For $g:\A \raa\B$,  $f_1:\B \raa\C_1$ and $f_2:\B \raa\C_2,$ we have
$$\la f_1,f_2 \ra \circ g \approx \la f_1 \circ g, f_2 \circ g \ra.$$
In terms of procedures, this says that doing $g$ and then doing both $f_1$ and $f_2$ is
the same as doing both $f_1 \circ g$ and $f_2 \circ g$, i.e., the following two flowcharts
are essentially the same.

$$\xymatrix{ &*+[F-,]{g} \ar[dl] \ar[dr] &&&\approx& *+[F-,]{g}\ar[d] &*+[F-,]{g}\ar[d]
\\
*+[F-,]{f_1}&&*+[F-,]{f_2}&&&*+[F-,]{f_1}&*+[F-,]{f_2}
}$$

In terms of trees, this amounts to saying that these trees are equivalent:

\Tree [ {$g:\A\raa\B$} [ {$f_1:\B\raa\C_1$} {$f_2:\B\raa\C_2$}
 ].{$\la f_1, f_2 \ra:\B\raa\C_1 \times \C_2$\\ \b} ].{$\la f_1, f_2 \ra \circ g:\A\raa\C_1 \times \C_2$\\ \c}
\Tree [ [{$g:\A\raa\B$} {$f_1:\B\raa\C_1$} ].{$f_1 \circ g:\A\raa\C_1$\\ \c}
[{$g:\A\raa\B$} {$f_2:\B\raa\C_2$} ].{$f_2 \circ g:\A\raa\C_2$\\ \c}
].{$\la f_1\circ g, f_2\circ g \ra:\A\raa\C_1 \times \C_2$\\ \b}

It is important to realize that it does not make sense to require composition to distribute
over bracket on the left:
$$g \circ \la f_1,f_2 \ra   \nsim \la g \circ f_1  , g \circ f_2  \ra.$$
The following two flowcharts are {\it not} essentially the same.
$$\xymatrix{*+[F-,]{f_1}\ar[dr] &&*+[F-,]{f_2}\ar[dl] &&\nsim& *+[F-,]{f_1}\ar[d] &*+[F-,]{f_2}\ar[d]
\\
&*+[F-,]{g}&&&&*+[F-,]{g}&*+[F-,]{g}
}$$
The left $g$ requires two inputs. The right $g$'s only require one.

\subsection{Bracket}
\noindent{\bf Bracket is Associative.}
The bracket is associative. For any three maps $f,g,$ and $h$ with the same domain, we have
$$ \la \la f,g \ra, h \ra \approx\la f,\la g, h \ra \ra $$
In terms of trees, this amounts to

\Tree [ [ {$f:\A\raa\B$} {$g:\A\raa\C$} ] .{$\la f,g \ra :\A\raa\B\times \C$ \\ \b} {$h:\A \raa\D$}
].{$\la \la f,g\ra h\ra:\A\raa\B\times \C\times \D$\\ \b}
$\approx$
\Tree [ {$f:\A\raa\B$} [ {$g:\A\raa\C$} {$h:\A \raa\D$} ].{$\la g,h \ra:\B\raa\C \times \D$\\ \b}
].{$\la f, \la g,h \ra \ra:\A\raa\B\times \C\times \D$\\ \b}

\vspace{.2in}

\noindent{\bf Bracket is Almost Commutative.}
It is not essential what is
written in the first or the second place. For any two maps $f$ and $g$ with the same domain,
$$\la f,g \ra \approx tw \circ \la g, f \ra.$$
In terms of trees, this amounts to

\Tree [{$f:\A\raa\B$} {$g:\A\raa\C$} ] .{$\la f,g \ra :\A\raa\B\times \C$ \\ \b}
$\approx$
\Tree [ [ {$g:\A\raa\C$} {$f:\A\raa\B$} ] .{$\la g,f \ra :\A\raa\C\times \B$ \\ \b} {$tw:\C \times \B\raa\B \times\C$}
].{$tw \circ \la g,f \ra: \A \raa\B \times \C$ \\ \c}

\vspace{.2in}

\noindent{\bf Twist is Idempotent.}
There are other relations that the twist map must respect.
Idempotent means
$$tw_{\A,\B} \circ tw_{\A,\B} \approx \pi^{\A \times \B}_{\A \times \B}:\A \times \B \raa\A \times \B.$$

\vspace{.2in}
\noindent{\bf Twist is Coherent.}
We would like the twist maps of three elements to get along with themselves.
$$(tw_{\B,\C}\times \pi_\A) \circ (\pi_\B \times tw_{\A, \C})\circ (tw_{\A,\B}\times \pi_\C)
\approx
(\pi_\C \times tw_{\A,\B})\circ (tw_{\A,\C}\times \pi_\B)\circ (\pi_\A \times tw_{\B,\C}).$$
This is called the hexagon law or the third Reidermeister move. Given the idempotence and hexagon laws,
it is a theorem
that there is a unique twist map made of smaller twist maps between any two products of elements (\cite{MacLane} Section
XI.4).

\vspace{.2in}

\noindent{\bf Bracket and Projections.}
A bracket followed by a projection onto the first output means the second output is ignored:
$ f \approx \pi^{\B \times \C}_\B \circ \la f,g \ra.$
In terms of trees, this amounts to

$f:\A\raa\B \qquad \approx$
\Tree [ [{$f:\A\raa\B$} {$g:\A\raa\C$} ] .{$\la f,g \ra :\A\raa\B\times \C$ \\ \b}
{$\pi^{\B\times\C}_\B:\B\times \C \raa\B$} ].{$\pi^{\B\times\C}_\B \circ \la f,g \ra :\A\raa\B$ \\ \c}


Similarly for a projection onto the second output:
$ g \approx \pi^{\B \times \C}_\C \circ \la f,g \ra.$
%
%

\vspace{.2in}
\noindent{\bf Bracket and Identity.} We want the bracket to be functorial, i.e., to respect the identity.
$$ \la \pi^\A_\A, \pi^\A_\A \ra \approx \triangle : \A \lra \A \times \A$$

\subsection{Bracket and Recursion}
When there are two unrelated processes, we can perform both of them in one loop or we can perform each of them
in its own loop.

\fbox{\parbox{1.9in}{\it  $h=\la f_1(x),f_2(x)\ra$ \\ For i = 1 to n \\ $~~~~h=(g_1(x,\pi_1 h), g_2(x,\pi_2 h ))$}}
$\approx$
\fbox{\fbox{\parbox{1.5in}{\it  $h_1=f_1(x)$ \\ For i = 1 to n \\ $~~~~h_1=g_1(x,h_1)$}};
      \fbox{\parbox{1.5in}{\it  $h_2=f_2(x)$ \\ For i = 1 to n \\ $~~~~h_2=g_2(x,h_2)$}}}

In $\sharp$ notation this amounts to saying
$$h=\la f_1,f_2 \ra \sharp (g_1 \boxtimes g_2) \approx \la f_1 \sharp g_1, f_2 \sharp g_2 \ra=\la h_1, h_2 \ra.$$

In terms of trees this says that this tree:

\Tree [ [{$f_1:\A\raa\B$} {$f_2:\A\raa\B$} ].{$\la f_1,f_2 \ra:\A \raa\B \times \B $\\ \b}
[ {$g_1:\A \times \B \raa\B $} {$g_2:\A \times \B \raa\B $}
] .{$g_1 \boxtimes g_2: \A \times \B \times \B \raa\B \times \B $\\ \boxt }
] .{$h=(\la f_1,f_2 \ra \sharp (g_1 \boxtimes g_2)):\A \times \N \raa\B \times \B$\\ \r }

\vspace{.2in}

\noindent is equivalent ($\approx $) to this tree:

\vspace{.2in}

\Tree [ [{$f_1:\A\raa\B$} {$g_1:\A \times \B \raa\B $}
].{$h_1=(f_1 \sharp g_1) :\A \times \N \raa\B$\\ \r}
[{$f_2:\A\raa\B$} {$g_2:\A \times \B \raa\B. $}
].{$h_2=(f_2 \sharp g_2) :\A \times \N \raa\B$\\ \r}
].{$\la h_1, h_2 \ra=\la f_1 \sharp g_1, f_2 \sharp g_2 \ra:\A \times \N \raa\B \times \B$\\ \b }

\subsection{Recursion and Composition}
\noindent{\bf Unwinding a Recursive Loop.}
Consider the following two algorithms

\fbox{\parbox{1.5in}{\it  $h=f(x)$ \\For i = 1 to n \\ $~~~~h=g_1(x,h)$ \\ $~~~~h=g_2(x,h)$ }} $\qquad$
\fbox{\parbox{1.5in}{\it $h'=g_1(x,f(x))$ \\ For i = 1 to n-1 \\ $~~~~h'=g_2(x,h')$ \\
$~~~~h'=g_1(x,h')$ \\ $h'=g_2(x,h')$ }}

This is the most general form of unwinding a loop. If $g_1$ is the identity process (does nothing), these become

\fbox{\parbox{1.5in}{\it $h=f(x)$\\For i = 1 to n \\   $~~~~h=g_2(x,h)$ }} $\qquad$
\fbox{\parbox{1.5in}{\it $h'=f(x)$\\For i = 1 to n-1 \\ $~~~~h'=g_2(x,h')$ \\ $h'=g_2(x,h').$ }}

If $g_2$ is the identity process, these become

\fbox{\parbox{1.5in}{\it $h=f(x)$\\For i = 1 to n \\  $~~~~h=g_1(x,h)$ }} $\qquad$
\fbox{\parbox{1.5in}{\it  $h'=g_1(x,f(x))$ \\ For i = 1 to n-1 \\   $~~~~h'=g_1(x,h').$  }}

In terms of recursion, the most general form of unwinding a loop, the left top box coincides with

$h(x,0) = f(x)$

$h(x,n+1) = g_2 (x,g_1(x,h(x,n))).$

\noindent The right top box coincides with:

$h'(x,0) = g_1(x,f(x))$

$h'(x,n+1) = g_1 (x,g_2(x,h'(x,n))).$

\noindent How are these two recursions related?
We claim that for all $n \in \N$ we have $g_1(x,h(x,n)) = h'(x,n).$
This may be proven by induction. The $n=0$ case is trivial. Assume it is true for $k$, and we shall show
it is true for $k+1$.
$$g_1(x,h(x,k+1))= g_1(x,g_2(x,g_1(x,h(x,k))))=g_1(x,g_2(x,h'(x,k)))=h'(x,k+1).$$
The first equality is from the definition of $h$; the second equality is the induction hypothesis; and the
third equality is from the definition of $h'$.

Although $g_1 \circt h$ and $g_2$ are constructed differently, they are essentially the
same program so we shall set them
equivalent to each other: $g_1 \circt h\approx h' $
If one leaves out the $h$ and $h'$ and uses the $\sharp$ notation, this becomes
$$g_1 \circt (f \sharp (g_2 \circt g_1)) \approx (g_1 \circt f)\sharp (g_1 \circt g_2).$$

In terms of trees, this means that

\Tree [ [{$f:\A \raa\B$} [{$g_2:\A \times \B \raa\B$} {$g_1:\A \times \B \raa\B$}
].{$g_2 \circt g_1: \A \times \B \raa\A$\\ \boxcirct} ].{$h:\A \times \N \raa\B$\\ \r} {$g_1:\A \times \B \raa\B$}
].{$g_1 \circt h: \A \times \N \raa\B$\\ \boxcirct}

\vspace{.2in}

\noindent is equivalent ($\approx$) to

\vspace{.2in}

\Tree [ [  {$f:\A \raa\B$} {$g_1:\A \times \B \raa\B$} ].{$g_1 \circt f:\A \raa\B $\\ \boxcirct}
[ {$g_2:\A \times \B \raa\B$} {$g_1:\A \times \B \raa\B$} ].{$g_1 \circt g_2:\A \times \B \raa\B$\\ \boxcirct}
].{$h':\A \times \N \raa\B $\\ \r}

\vspace{.2in}

\noindent{\bf Recursion and Null.}
If $h$ is defined by recursion from $f$ and $g$, i.e. $h=f \sharp g$,
then by definition of recursion
$h(x,0)=f(x)$
or
$h(x,n(y))=f(x)$
where $n$ is the null function and $y \in \N$. This means $h \circt n = f$.
We shall set these equivalent
$h \circt n \approx f$
Using the $\sharp$ notation, this amounts to:
$(f \sharp g) \circt n \approx f.$
In terms of algorithms, this amounts to saying that the following two algorithms are
equivalent:

\fbox{\parbox{1.5in}{\it  $h=f(x)$\\ For i= 0 to 0 \\ $~~~~h=g(x,h)$ }}
$\approx$ \fbox{\parbox{1.5in}{\it  $h=f(x)$ }}

In terms of trees, this is

\Tree [ {$n:\N \raa\N$} [{$f:\A \raa\B$} {$g:\A \times \B \raa\B$}
].{$h:\A \times \N \raa\B$\\ \r} ].{$(h \circt n):\A \raa\B$\\ \boxcirct}
$\approx \qquad $
$f:\A \raa\B$

Notice that the $g$ on the left tree is not on the right tree.

\vspace{.2in}
\noindent{\bf Recursion and Successor.}
Let $h$ be defined by recursion from $f$ and $g$, i.e., $h=f \sharp g$.
Then by definition of recursion:
$h(x,k+1)=g(x,h(x,k))$
or
$h(x,s(k))=g(x,h(x,k))$
where $s$ is the successor function and $k \in \N$. This is the same as $h \circt s = g \circt h$.
We shall set them equivalent
$h \circt s \approx g \circt h.$
Using the $\sharp$ notation, this becomes $(f \sharp g) \circt s \approx g \circt (f \sharp g). $
In terms of algorithms, this says that the following two algorithms are equivalent

\fbox{\parbox{1.5in}{\it $h=f(x)$ \\ For i = 1 to k+1 \\ $~~~~h=g(x,h)$ }} $\approx$
\fbox{\parbox{1.5in}{\it $h=f(x)$ \\ For i = 1 to k \\ $~~~~h=g(x,h)$\\
$h=g(x,h)$ }}

\noindent In terms of trees, this says that the following two trees are set equivalent

\Tree [ {$s:\N \raa\N$} [ {$f:\A \raa\B$} {$g:\A \times \B \raa\B$} ].{$h:\A \times \N \raa\B$\\ \r}
].{$h \circt s:\A \times \N \raa\N $\\ \boxcirct}
$\approx$
\Tree [ [ {$f:\A \raa\B$} {$g:\A \times \B \raa\B$} ].{$h:\A \times \N \raa\B$\\ \r}
{$g:\A \times \B \raa\B$} ].{$g \circt h:\A \times \N \raa\N $ \\ \boxcirct}

\vspace{.2in}
\noindent{\bf Recursion and Identity.} If $g=\pi^{\A \times \B}_\B$, i.e., if we do recursion over the
identity function, then we are not really doing recursion at all.
$$(f \sharp g) = (f\sharp \pi^{\A\times \B}_\B) \quad \approx \quad (f \circ \pi^{\A \times \N}_\A).$$

\subsection{Products}
\noindent{\bf The product is associative.} That is for any three maps $f:\A \raa\A',$ $g:\B \raa\B'$ and
$h:\C \raa\C'$ the two products are equivalent:
$$f \times (g \times h) \approx (f \times g) \times h: \A\times \B \times \C \raa\A' \times \B' \times \C'.$$
This follows immediately from the associativity of bracket.

\vspace{.2in}
\noindent{\bf The product respects identity.}
$$ \pi^\A_\A \times \pi^\B_\B \approx  \pi^{\A \times \B}_{\A \times \B}. $$
This falls out of the fact that the bracket respects the identity.

\vspace{.2in}
\noindent {\bf Interchange Rule.} We must show that the product and the composition respect
each other. In terms of maps, this corresponds to the following situation:

$$\xymatrix{
\A_1 \ar[dd]_{f_1} \ar@/_.5in/[dddd]_{f_2 \circ f_1}&&
\A_1 \times B_1 \ar[ll]_{\pi} \ar[rr]^{\pi}\ar[dd]_{f_1 \times g_1}&& \B_1\ar[dd]_{g_1}\ar@/^.5in/[dddd]^{g_2 \circ g_1}
\\
\\
\A_2 \ar[dd]_{f_2}&& \A_2 \times B_2 \ar[ll]_{\pi} \ar[rr]^{\pi}\ar[dd]_{f_2 \times g_2}&& \B_2\ar[dd]_{g_2}
\\
\\
\A_3 && \A_3 \times B_3 \ar[ll]_{\pi} \ar[rr]^{\pi}&& \B_3
}$$

$(f_2 \times g_2) \circ (f_1 \times g_1)$ and $(f_2 \circ f_1) \times (g_2 \circ g_1)$
are two ways of getting from $\A_1 \times \B_1$ to $\A_3 \times \B_3$. We shall declare
these two methods equivalent:
$$(f_2 \times g_2) \circ (f_1 \times g_1) \approx (f_2 \circ f_1) \times (g_2 \circ g_1).$$

In terms of trees, this tree:

\Tree [ [{$f_1:\A_1 \raa\A_2 $} {$g_1:\B_1 \raa\B_2 $}
].{$f_1 \times g_1:\A_1 \times \B_1 \raa\A_2 \times \B_2$ \\ \p}
[{$f_2:\A_2 \raa\A_3 $} {$g_2:\B_2 \raa\B_3$}
 ].{$f_2 \times g_2:\A_2 \times \B_2 \raa\A_3 \times \B_3 $ \\ \p}
].{$(f_2 \times g_2) \circ (f_1 \times g_1) :\A_1 \times \B_1 \raa\A_3 \times \B_3$\\ \c}

\vspace{.2in}

\noindent is equivalent ($\approx$) to this tree:

\vspace{.2in}

\Tree [ [{$f_1:\A_1 \raa\A_2 $} {$f_2:\A_2 \raa\A_3 $} ].{$f_2 \circ f_1:\A_1 \raa\A_3 $\\ \c}
[ {$g_1 :\B_1 \raa\B_2$} {$g_2 :\B_2 \raa\B_3.$} ].{$g_2 \circ g_1:\B_1 \raa\B_3 $\\ \c}
].{$(f_2 \circ f_1) \times (g_2 \circ g_1):\A_1 \times \B_1 \raa\A_3 \times \B_3 $\\ \p}

\vspace{.2in}

One should realize that this equivalence is not anything new added to our list of equivalences. It is actually
a consequence of the definition of product and the equivalences that we assume about bracket. In detail
$$(f_2 \times g_2) \circ (f_1 \times g_1)=\la f_2 \pi, g_2 \pi \ra \circ \la f_1 \pi, g_1 \pi \ra \approx
\la f_2 \pi \la f_1 \pi , g_1 \pi \ra \ra, g_2 \pi \la f_1 \pi, g_1 \pi\ra \ra$$
$$\approx \la f_2 \circ f_1 \pi, g_2 \circ g_1 \pi \ra =  (f_2 \circ f_1) \times (g_2 \circ g_1).$$
The first and the last equality are from the definition of product. The first equivalence comes from the
fact that composition distributes over bracket. The second equivalence is a consequence of the relationship
between the projection maps and the bracket.

\section{Algorithms}
We have given relations telling when two programs/trees/descriptions are similar. We would like to look at the equivalence
classes that these relations generate. It will become apparent that by taking $\PRdesc$ and ``modding out'' by these equivalence relations,
we shall get more structure.

The relations split up into two disjoint sets: those for which there is a loss of information and those for
which there is no loss of information. Let us call the former set of relations {\bf (I)} and the latter set
{\bf (II)}. The following relations are in group {\bf (I)}.
\begin{enumerate}
\item Null Function and Composition: $n \circ f \approx n \circ \pi^\A_\N$
\item Bracket and First Projection: $ f \approx \pi^{\B \times \C}_\B \la f,g \ra$
\item Bracket and Second Projection: $ g \approx \pi^{\B \times \C}_\C \la f,g \ra$
\item Recursion and Null Function: $(f \sharp g) \circt n \approx f$
\end{enumerate}

After setting these trees equivalent, there exists the following quotient graph and graph morphism.
$$\xymatrix{\PRdesc \ar@{>>}[rr] && \PRdesc /({\bf I})  }$$
In detail, $\PRdesc /({\bf I})$ has the same vertices as $\PRdesc$, namely powers of the set of natural numbers. The
edges are equivalence classes of edges of $\PRdesc$.

Descriptions of primitive recursive functions which are equivalent to ``pruned'' descriptions
by relations of type {\bf (I)} we shall call ``stupid descriptions''. They are descriptions that are wasteful in
the sense that
part of their tree is dedicated to describing a certain function and that function is not needed. The part of the
tree that describes the unneeded function can be lopped off.
One might call $\PRdesc /({\bf I})$
the graph of ``intelligent descriptions'' since within this graph every ``stupid descriptions'' is equivalent
to another program without the wastefulness.

We can further quotient
$\PRdesc /({\bf I})$ by relations of type {\bf (II)}:

\begin{enumerate}
\item Composition Is Associative: $f \circ (g \circ h) \approx (f \circ g) \circ h.$
\item Projections Are Identities: $f \circ \pi^\A_\A \approx f \approx \pi^\B_\B \circ f.$
\item Composition Distributes Over Bracket: $\la f_1,f_2 \ra \circ g \approx \la f_1 \circ g, f_2 \circ g \ra.$
\item Bracket Is Associative: $\la f,\la g, h \ra \ra \approx \la \la f,g \ra, h \ra .$
\item Bracket Is Almost Commutative: $\la f,g \ra \approx tw \circ \la g, f \ra.$
\item Bracket is functorial: $ $$ \la \pi^\A_\A, \pi^\A_\A \ra \approx \triangle$$
 $
\item Twist Is Idempotent: $tw \circ tw = \pi.$
\item Reidermeister III: $$(tw_{\B,\C}\times \pi_\A) \circ (\pi_\B \times tw_{\A, \C})\circ (tw_{\A,\B}\times \pi_\C)
\approx
(\pi_\C \times tw_{\A,\B})\circ (tw_{\A,\C}\times \pi_\B)\circ (\pi_\A \times tw_{\B,\C}).$$
\item Recursion and Bracket: $\la f_1,f_2 \ra \sharp (g_1 \boxtimes g_2) \approx \la f_1 \sharp g_1, f_2 \sharp g_2 \ra.$
\item Recursion and Composition: $g_1 \circt (f \sharp (g_2 \circt g_1)) \approx (g_1 \circt f)\sharp (g_1 \circt g_2).$
\item Recursion and Successor Function: $(f \sharp g) \circt s \approx g \circt (f \sharp g). $
\end{enumerate}

There is a further projection onto the quotient graph:
$$\xymatrix{\PRdesc \ar@{>>}[r] & \PRdesc /({\bf I}) \ar@{>>}[r] & \PRalg
= (\PRdesc/{\bf I})/{\bf II} =\PRdesc / (({\bf I})\bigcup({\bf II})).  }$$
$\PRalg$, or primitive recursive algorithms, are the main object of interest in this Section.

What does $\PRalg$ look like? Again the objects are the same as $\PRdesc$, namely powers of the set of natural numbers. The
edges are equivalence classes of edges of $\PRdesc$.

What type of structure does it have? In $\PRalg$, for any three composable arrows, we have
$$f \circ (g \circ h) = (f \circ g) \circ h$$
and for any arrow $f:\A \raa\B$ we have
$$f \circ \pi^\A_\A = f = \pi^\B_\B \circ f.$$
That means that composition is associative and that the $\pi$'s act as identities. Whereas $\PRdesc$ was
only a graph with a composition and identities that did not act like identities, $\PRalg$ is a genuine category.

$\PRalg$ has more structure than only a category. For one, there is a strictly associative product. On objects, the product structure is obvious:
$$\N^m \times \N^n = \N^{m+n}.$$
On morphisms, the product $\times$ was defined using the bracket above. The $\pi$ are the projections of the product.
In $\PRalg$ the twist map is idempotent and coherent.
The fact that the product respects the composition is expressed with the interchange rule.

The category $\PRalg$ is closed under recursion. In other words, for
any $f:\A \raa\B$ and any $g:\A \times \B \raa\B$, there exists an
$h:\A \times \N \raa\B$ defined by recursion. The categorical way of
saying that a category is closed under recursion, is to say that the
category contains a weak parameterized natural number object. The
simplest definition of a weak natural number object in a category
is a diagram
$$\xymatrix{\ast \ar[rr]^0  && \N \ar[rr]^s  && \N }$$
such that for any $k\in \N$ and $g:\N \raa\N$, there exists an $h:\N \raa\N$ such that the
following diagram commutes.
$$\xymatrix{
\ast \ar[rr]^0 \ar[rrdd]_k && \N \ar[rr]^s \ar[dd]_h && \N \ar[dd]^h \\
\\
&& \N \ar[rr]_g && \N }$$ (See e.g. \cite{TTT1,TTT,MacLane}).
Following \cite{Lambekandscott}, we do not insist that the $h$ is a
unique morphism that satisfies the condition. When there is
uniqueness, we say that the natural number object is strong. Saying
that the above diagram commutes is the same as saying that $h$ is
defined by the simplest recursion scheme. For our more general
version of recursion, we require a weak {\it parameterized} natural
number object, that is, for every $f:\A \raa\B$ and $g:\A \times
\B\raa\B$ there exists a $h: \A \times \N \raa\B$ such that the
following two squares commute.
$$\xymatrix{
\A \times \ast \ar[rr]^{\pi \times 0} \ar[dd]_\wr && \A \times \N \ar[dd]^h\\
\\
\A \ar[rr]_f && \B
}
\qquad
\xymatrix{
\A \times \N \ar[rr]^{\pi \times s} \ar[dd]_{\la \pi^{\A \times \N}_{\A},h\ra} && \A \times \N \ar[dd]^h\\
\\
\A \times \B \ar[rr]_g && \B
}$$

From the fact that in $\PRalg$ we have an object $\N$, the morphisms
$0:\ast \raa\N$ and $s:\N \raa\N$ and these morphisms satisfy
$h \circt n = (f \sharp g) \circt n = f$ and $h \circt s = (f \sharp g) \circt s = g \circt (f \sharp g)=g \circt h,$
we see that $\PRalg$ has a weak parameterized natural number object.

Some words on the uniqueness of $h$ are needed. Given descriptions
$f$ and $g$ of the correct arity, we can form the description $h=(f
\sharp g)$. This $h$ will satisfy the requirements of the
parameterized natural number object. But there is no reason to think
that this is the only description that would satisfy the
requirements. Any other description of the same function that $h$
performs would also satisfy the requirement. This is in sharp
contrast to a category of functions. Given primitive recursive
functions $f$ and $g$ of the right arity, there is {\it only one}
function $h=(f\sharp g)$ that satisfies the recursion axiom. One can
think of this distinction as a fundamental difference between syntax and
semantics. In a syntactical category, it is impossible to demand
uniqueness. There are many descriptions of objects that satisfy
conditions. In contrast, within semantic categories, there is only
one object that satisfies requirements. In Lambek and Scott
\cite{Lambekandscott}, they deal with syntactical categories of
proof and there too, they only have a weak natural number objects
(page 46). Similarly, in Peter Johnstone's discussion of
lambda-calculus in Proposition 4.2.12 on page 959 of volume II of
\cite{Johnstone}, the natural number object in the syntactical
category is weak.

We must show that in $\PRalg$, the natural number object respects the bracket operation.
This fundamentally says that the central square in the following two diagrams commute.
$$\xymatrix{
&&\A \times \ast \ar[rrrr]^{\pi \times 0} \ar[ddr]_\wr \ar[dddd]_\wr \ar[ddddddll]_\wr &&&&
\A \times \N
\ar[ddl]^{h_1} \ar[dddd]^{\la h_1,h_2 \ra} \ar[ddddddrr]^{h_2}
\\
\\
&&&\A \ar[rr]_{f_1} && \B
\\
\\
&&\A \ar[rrrr]_{\la f_1,f_2\ra} \ar@{=>}[ddll] \ar@{=>}[uur]&&&& \B \times \B \ar[uul]^\pi \ar[ddrr]_\pi
\\
\\
\A \ar[rrrrrrrr]_{f_2} &&&&&&&& \B
}$$
The left hand triangles commute from the fact that $\ast$ is a terminal object. The right hand triangles commute
because the equivalence relation forced the projections to respect the bracket. The inner and
outer quadrilateral are assumed to commute. We conclude that the central square commutes.
$$\xymatrix{
&&\A \times \N \ar[rrrr]^{\pi \times s} \ar[ddr]_\wr \ar[dddd]_\wr \ar[ddddddll]_\wr &&&&
\A \times \N
\ar[ddl]^{h_1} \ar[dddd]^{\la h_1,h_2 \ra} \ar[ddddddrr]^{h_2}
\\
\\
&&&\A \times \B  \ar[rr]_{g_1} && \B
\\
\\
&&\A \times \B \times \B \ar[rrrr]_{g_1 \boxtimes g_2}\ar[uur]^\pi
\ar[ddll]^\pi &&&& \B\times \B \ar[uul]^\pi \ar[ddrr]_\pi
\\
\\
\A \times \B \ar[rrrrrrrr]_{g_2} &&&&&&&& \B
}$$
Similarly, the left and the right triangles commute because the projections act as they are supposed to.
The inner and outer quadrilateral commute out of assumption. We conclude that central square commutes.

We also must show that the natural number object respects the composition of morphisms.
In $\sharp$ notation this amounts to
$$g_1 \circt (f \sharp (g_2 \circt g_1)) = (g_1 \circt f)\sharp (g_1 \circt g_2).$$
For the simpler form of recursion, this reduces to
$$g_1 \circ (k \sharp (g_2 \circ g_1)) = (g_1 \circ k)\sharp (g_1 \circ g_2).$$
Setting $h=k \sharp (g_2 \circ g_1)$ and $h'=(g_1 \circ k)\sharp (g_1 \circ g_2)$, we get
the following natural number object diagram
$$\xymatrix{
\ast \ar[rr]^0 \ar[rddd]^k && \N \ar[rrrr]^s \ar[dddl]_h \ar[dddr]_{h'}&&&&\N \ar[dddl]_h \ar[dddr]_{h'}
\\
\\
\\
&\B \ar[rr]_{g_1} &&\B \ar[rr]_{g_2}&& \B \ar[rr]_{g_1}&& \B.}$$
With the properties of $h$ and $h'$ we get that the triangles
commute.

Once we have $\PRalg$, we might ask when do two algorithms perform the same operation. We make
an equivalence relation and say two algorithms are equivalent ($\approx'$) iff they perform
the same operation. By taking a further quotient of $\PRalg$ we get $\PRfunc$. What does $\PRfunc$
look like. The objects are again powers of the set of natural numbers and the morphisms are
primitive recursive functions.

In summary, we have the following diagram.
$$\xymatrix{\PRdesc \ar@{>>}[r]& \PRdesc /({\bf I}) \ar@{>>}[r]&
\PRalg = \PRdesc / (({\bf I})\bigcup ({\bf II})) \ar@{>>}[r]
&\PRfunc = \PRalg / \approx'.}$$

Let us spend a few moments with some category theory. There is the
category $\cat$ of all (small) categories and functors between them.
Consider also the category $\catxn$. The objects are triples, $(\C,
\times, N)$ where $\C$ is a (small) category, $\times$ is a strict
product on $\C$ and $N$ is a weak parameterized natural number
object in $\C$. The morphisms of $\catxn$ are functors $F: (\C,
\times, N) \raa(\C', \times', N')$ that respect the product and
the natural number object. For $F: \C \raa\C'$ to respect the product,
we  mean that
$$\mbox{For all }f,g\in \C\quad F(f\times g) = F(f) \times' F(g).$$
To say that $F$ respects the natural number object means that
if $$\xymatrix{\ast \ar[rr]^0  && N \ar[rr]^s  && N }$$ is a natural number object in $\C$
and $$\xymatrix{\ast' \ar[rr]^{0'}  && N' \ar[rr]^{s'}  && N' }$$ is a natural number
object in $\C'$ then $F(N)=N', F(\ast)=\ast', F(0)=0'$ and $F(s)=s'$.
For a given natural number object in a category, there is an implied function $\sharp$ that takes
two morphisms $f$ and $g$ of the appropriate arity and outputs the unique $h=f\sharp g$ of the appropriate
arity. Our definition of a morphism between two objects in $\catxn$ implies that
$$\mbox{For all appropriate }f,g\in \C\quad F(f\sharp g) = F(f) \sharp' F(g).$$

There is an obvious forgetful functor $U:\catxn \raa\cat$ that takes $(\C, \times, N)$ to $\C$.
There exists a left adjoint to this forgetful functor:
$$\xymatrix{ \cat \ar@/^/[rr]^L_\bot && \catxn. \ar@/^/[ll]^U}$$
This adjunction means that for all small categories $\C \in \cat$ and $\D\in \catxn$
there is an isomorphism
$$\catxn (L(\C), \D) \simeq \cat (\C, U(\D)).$$
Taking $\C$ to be the empty category ${\bf \emptyset}$ we have
$$\catxn (L({\bf \emptyset}), \D) \simeq \cat ({\bf \emptyset}, U(\D)).$$
Since ${\bf \emptyset}$ is the initial object in $\cat$, the right
set has only one object. In other words $L({\bf \emptyset})$ is a
free category with product and a weak parameterized natural number
object and it is an initial object in the category $\catxn$.

We claim that $L({\bf \emptyset})$ is none other then our category $\PRalg$.

\begin{teo}
$\PRalg$ is an initial object in the category of categories with a
strict product and a weak parameterized natural number object.
\end{teo}

We have already shown that $\PRalg$ is a category with a strict product and a natural number object. It remains to
be shown that for any object $(\D, \times, N')\in \catxn$ there is a unique functor
$F_\D:\PRalg \raa\D$. Our task is already done by recalling that the objects and morphisms in
$\PRalg$ are all generated by the natural number object and that functors in $\catxn$ must preserve
this structure. In detail, $F_\D(\N)=N'$ and since $F_\D$ must preserve products $F_\D(\N^i)=(N')^i$.
And similarly for the morphisms of $\PRalg$. The morphisms are generated by the $\pi$s, the $n$ and $s$ in
the natural number object of $\PRalg$. They are generated by composition, product and recursion. $F_\D$ is a functor
and so it preserves composition. We furthermore assume it preserves product and recursion. $(\D, \times, N')\in \catxn$
might have many more objects and morphisms but that is not our concern here. $\PRalg$ has very few morphisms.

The point of this theorem is that $\PRalg$ is not simply a nice category where all algorithms live. Rather it is
a category with much structure. The structure tells us how algorithms are built out of each other.
$\PRalg$ by itself is not very interesting. It is only its extra structure that demonstrates the importance
of this theorem. $\PRalg$ is not simply the category made of algorithms, rather, it is the category that
makes up algorithms.

$\PRfunc$ is the smallest category with a strict product and a strong
parameterized natural number object.

Before we go on to other topics, it might be helpful to ---literally--- step away from the trees and look
at the entire forest. What did we do here? The graph $\PRdesc$ has operations. Given edges of the
appropriate arity, we can compose them, bracket them or do recursion on them. But these operations
do not have much structure. $\PRdesc$ is not even a category. By placing equivalence relations on $\PRdesc$,
which are basically coherence relations, we are giving the quotient category better and more amenable
structure. So coherence theory, sometimes called higher-dimensional algebra, tells us when
two programs are essentially the same.

\section{Complexity Results}
An algorithm is not {\it one} arrow in the category $\PRalg$. An algorithm is a scheme of arrows,
one for every input size. We need a way of choosing each of these arrows.

There are many different
species of algorithms. There are algorithms that accept $n$ numbers and output one number. A scheme for
such an algorithm might look like this:
$$\xymatrix{
&& \N^1\ar[dd]^{c_1}\\
&&& \N^2 \ar[dl]^{c_2}\\
\vdots \ar[rr]^{c}&& \N && \N^3\ar[ll]^{c_3}\\
&\N^k \ar[ur]^{c_k}&& \N^4\ar[ul]^{c_4}\\
&&\cdots\ar[uu]^{c}
}$$
We shall call such a graph a {\it star graph} and denote it $\bigstar$.

However there are other species of algorithms. There are algorithms that accept $n$ numbers and
output $n$ numbers (like sorting or reversing a list, etc.) Such a scheme looks like
$$\xymatrix{
\N^1 \ar[r]^{c_1}& \N^1&
\N^2 \ar[r]^{c_2}& \N^2&
 \dots &
\N^k \ar[r]^{c_k}& \N^k
& \dots &
}$$
We shall also call such a graph a star graph.

One can think of many other possibilities. For example, algorithms that accept $n$ numbers
and outputs their max, average and minimum (or mean, median and mode) outputs three numbers.
We shall not be particular as to what what type of star graph we will be working with.

Given any star graph $\bigstar$, a scheme that chooses one primitive
recursive description for each edge is a graph homomorphism
$Sch:\bigstar \raa\PRdesc$ that is the identity on vertices, i.e.,
$Sch(\N^i)=\N^i$ for all $i \in \N$.

Composing $Sch:\bigstar \raa\PRdesc$ with the projection onto the equivalence classes $\PRdesc \raa\PRdesc/{\bf (I)}$
gives a graph homomorphism $\bigstar \raa\PRdesc/{\bf (I)}$. In order not to have too many names flying around, we shall
also call this graph homomorphism $Sch$. Continuing to compose with the projections, we get the following
commutative diagram.
$$\xymatrix{ &\bigstar\ar[dr]^{Sch} \ar[d]^{Sch}  \ar[dl]^{Sch}  \ar[drr]^{Sch}
\\\PRdesc \ar@{>>}[r]&\PRdesc /({\bf I}) \ar@{>>}[r]
  &\PRalg \ar@{>>}[r]
& \PRfunc. }$$

We are not interested in only one graph homomorphism $\bigstar \raa\PRdesc$. Rather we are interested
in the set of all graph homomorphisms. We shall call this set $\PRdesc^\bigstar$. Similarly, we shall look
at the set of all graph homomorphisms from $\bigstar$ to $\PRdesc/{\bf (I)}$, which we shall denote
$(\PRdesc/{\bf (I)})^\bigstar$. There is also $\PRalg^\bigstar$ and $\PRfunc^\bigstar$. There are also
obvious projections:
$$\xymatrix{\PRdesc^\bigstar \ar@{>>}[r] & (\PRdesc /({\bf I}))^\bigstar  \ar@{>>}[r]
&
\PRalg^\bigstar  \ar@{>>}[r]
& \PRfunc^\bigstar  }$$

\vspace{.5in}

Perhaps it is time to get down from the abstract highland and give two examples. We shall present mergesort and
insertion sort as primitive recursive algorithms. They are two different members of
$\PRalg^\bigstar$. These two different algorithms perform the same function in $\PRfunc^\bigstar$.

\vspace{.2in}

\begin{examp}
Mergesort depends on an algorithm that merges two sorted lists into one sorted list. We define an algorithm
$Merge$ that accepts $m$ numbers of the first list and $n$ numbers of the second list. $Merge$ inputs and
outputs $m+n$ numbers.
$$Merge_{0,1}(x_1)=Merge_{1,0}(x_1)=\pi^1_1(x_1)=x_1$$
$$Merge_{m,n}(x_1,x_2,\ldots,x_m,x_{m+1},\ldots, x_{m+n})=$$
$$\left\{ \begin{array}{l@{\quad : \quad}l}
(Merge_{m,n-1}(x_1,x_2,\ldots,x_m,x_{m+1},\ldots, x_{m+n-1}),x_n) & x_m \leq x_n \\
(Merge_{m-1,n}(x_1,x_2,\ldots,x_{m-1},x_{m+1},\ldots, x_{m+n}),x_m) & x_m > x_n  \end{array} \right. $$
With $Merge$ defined, we go on to define $MergeSort$. $MergeSort$ recursively splits the list into two parts, sorts
each part and then merges them.
$$MergeSort_1(x)=\pi^\N_\N(x)=x$$
$$MergeSort_k(x_1,x_2,\ldots, x_k)=$$
$$Merge_{\llcorner k/2 \lrcorner, \ulcorner k/2 \urcorner}
(MergeSort_{\llcorner k/2 \lrcorner}(x_1,x_2,\ldots, x_{\llcorner k/2 \lrcorner}),
MergeSort_{\ulcorner k/2 \urcorner}(x_{\llcorner k/2 \lrcorner+1},x_{\llcorner k/2 \lrcorner +2}, \ldots, x_k)$$
We might write this in short as
$$MergeSort= Merge \circ \la MergeSort,  MergeSort \ra $$
\end{examp}

\vspace{.2in}

\begin{examp}
Insertion sort uses an algorithm $Insert:\N^k \times  \N \raa\N^{k+1}$ which takes an ordered list of $k$
numbers adds a $k+1$th number to that list in its correct position. In detail,
$$Insert_0(x)=\pi^1_1(x)=x$$
$$Insert_k(x_1,x_2,\ldots, x_k, x)=$$
$$\left\{ \begin{array}{l@{\quad : \quad}l}
(x_1,x_2,\ldots,x_k,x) & x_k \leq x \\
(Insert_{k-1}(x_1,x_2,\ldots,x_{k-1},x),x_k)
& x_k > x  \end{array} \right. $$
The top case is the function $\pi^k_k \times \pi^1_1$ and the bottom case is the function
$(Insert_{k-1}\times \pi) \circ (\pi^{k-1}_{k-1} \times tw_{\N,\N})$. With $Insert$ defined, we go on to
define $InsertionSort$.
$$InsertionSort_1(x)=\pi^\N_\N(x)=x$$
$$InsertionSort_k(x_1,x_2,\ldots, x_k) = Insert_{k-1}(InsertionSort_{k-1}(x_1,x_2,\ldots, x_{k-1}), x_k)$$
We might write this in short as
$$InsertionSort=Insert(InsertionSort \times \pi)$$
\end{examp}

The point of the these examples, is to show that although these two algorithms perform the same function, they
are clearly very different algorithms. Therefore one can not say that they are ``essentially'' the same.

\vspace{.5 in}

Now that we have placed  the objects of study in order, let us classify them via complexity theory.
The only operations in our trees that are of any complexity is the recursions. Furthermore, the recursions
are only interesting if they are nested within each other. So for a given tree that represents a description
of a primitive recursive function, we might ask what is the largest number of nested recursions in this tree.
In other words, we are interested in the largest number of ``R'' labels on a path
from the root to a leaf of the tree. Let us call this the $Rdepth$ of the tree.

Formally, $Rdepth$ is defined recursively on the set of our labeled
binary trees. The $Rdepth$ of a one element tree is $0$. The
$Rdepth$ of an arbitrary tree $T$ is
$$Rdepth(T)=Max\left\{ Rdepth(left(T)),Rdepth(right(T))\right\} + (label(T) == \r)$$
where $(label(T) == \r)=1$ if the label of the root of $T$ is $\r$, otherwise it is $0$.

It is known that a primitive recursive function that can be expressed by a tree with $Rdepth$
of $n$ or less is an element of Grzegorczyk's hierarchy class $\mathcal{E}^{n+1}$. (See \cite{Clote}, Theorem 3.31
for sources.)

Complexity theory deals with the partial order of all functions
$\{f|f:\N \raa\R^{+}\}$ where
$$f \leq g \mbox{ iff } Lim_{n \rightarrow \infty} \frac{f(n)}{g(n)}<\infty.$$

For every algorithm we can associate a function that describes the $Rdepth$ of the trees used in that
algorithm. Formally, for a given algorithm, $A:\bigstar \raa\PRdesc$, we can associate a function $f_A:\N \raa\R^{+}$
where
$$f_A(n)=Rdepth(A(c_n))$$ when $c_n$ is an edge in $\bigstar$.
The function $\PRdesc^\bigstar \raa\{f|f:\N \raa\R^{+}\}$ where
$A \mapsto f_A$ shall be called $Rdepth_0$.

We may extend $Rdepth_0$ to $$Rdepth_1:(\PRdesc/{\bf (I)})^\bigstar \raa\{f|f:\N \raa\R^{+}\}.$$
For a scheme of algorithms $[A]:\bigstar \raa(\PRdesc/{\bf (I)})$ we define
$$f_{[A]}(n)=Min_{A'} \{Rdepth(A'(c_n))\}$$
where the minimization is over all descriptions $A'$ in the equivalence class $[A]$. (For the categorical
cognoscenti, $Rdepth_1$ is a right Kan extension of $Rdepth_0$ along the projection $\PRdesc^\bigstar \lra
(\PRdesc/{\bf (I)})^\bigstar$.

$Rdepth_1$ can easily be extended to $$Rdepth_2:\PRalg^\bigstar \raa\{f|f:\N \raa\R^{+}\}.$$ The following
theorem will show us that we do not have to take a minimum over an entire equivalence class.

\begin{teo}
Equivalence relations of type {\bf (II)} respect $Rdepth$.
\end{teo}
\noindent{\bf Proof.} Examine all the trees that express these relations throughout this paper. Notice that
if two trees are equivalent, then  their $Rdepth$s are equal. $\Box$

$Rdepth_2$ can be extended to $$Rdepth_3:\PRfunc^\bigstar \raa\{f|f:\N \raa\R^{+}\}.$$
We do this again with a minimization over the entire equivalence class (i.e. a Kan extension.)

And so we have the following (not necessarily commutative) diagram.

$$\xymatrix{\PRdesc^\bigstar \ar[dr]_{Rdepth_0} \ar@{>>}[r]&
(\PRdesc /({\bf I}))^\bigstar \ar@{>>}[r] \ar[d]_{Rdepth_1}&
\PRalg^\bigstar \ar@{>>}[r] \ar[dl]_{Rdepth_2}&
\PRfunc^\bigstar \ar[dll]^{Rdepth_3} \\
&\{f|f:\N \raa\R^{+}\} }$$

\begin{cor}
The center triangle of the above diagram commutes.
\end{cor}
This is in contrast to the other two triangles which do not commute.

In order to see why the right triangle does not commute, consider an inefficient sorting algorithm. $Rdepth_2$
will take this inefficient algorithm to a large function $\N \raa\R^+$. However, there are efficient sorting
algorithms and $Rdepth_3$ will associate a smaller function to the primitive recursive function of sorting.

There are many subclasses of $\{f|f:\N \raa\R^{+}\}$ like polynomials or exponential functions. Complexity
theory studies the preimage of these subclasses under the function $Rdepth_3$. The partial order in
$\{f|f:\N \raa\R^{+}\}$ induces a partial order of subclasses of $\PRfunc$ which are the ``complexity classes.''

\section{Future Directions}

We are in no way finished with this work and there are many directions that it can be extended.

\noindent {\bf Extend to all Computable Functions.} The most obvious project that we are pursuing is to extend
this work from primitive recursive functions to
all computable functions. In order to do this we must add the minimization operation. For a given
$g:\A \times \N \raa\N$, there is an $h:\A \raa\N$ such that
$$h(x)=Min_n \left\{ g(x,n)=1 \right\}$$

Categorically, this amounts to looking at the total order of $\N$. This induces an
order on the set of all functions from $\A$ to $\N$. We then look at all functions $h'$
that make this square commute.
$$\xymatrix{
\A \ar[rr]^{!} \ar[dd]_{\la \pi^\A_\A ,h' \ra}&& \ast \ar[dd]^{1}
\\
\\
\A \times \N \ar[rr]^g && \N
}$$
i.e.,
$$g(x,h'(x))=1.$$
Let $h:\A \raa\N$ be the minimum such function.

We might want to generalize this operation. Let $f:\A \raa\B$ and
$g:\A \times \N \raa\B$, then we define $h:\A \raa\N$ to be the function
$$h(x)=Min_n \left\{ g(x,n)=f(x) \right\}.$$
Categorically, this amounts to looking at all functions $h'$ that make the triangle commute:
$$\xymatrix{
&\A \ar[ddr]^f \ar[ddl]_{\la \pi^\A_\A, h' \ra}
\\
\\
\A \times \N \ar[rr]_g &&\B
}$$
i.e.,
$$g(x,h'(x))=f(x).$$

Let $h:\A \raa\N$ be the minimum such function.

Hence minimization is a fourth fundamental operation:
\qtreecentertrue

\Tree [ {$f:\A \raa\B$} {$g:\A \times \N \raa\B$} ].{$h:\A \raa\N $\\ \m}

There are several problems that are hard to deal with. First, we leave the domain of total
functions and go into the troublesome area of partial functions. All the relational axioms
have to be reevaluated from this point of view. Second, what should we substitute for $Rdepth$
as a complexity measure?

Progress is being made in this direction in a forthcoming paper by
Yuri Manin and the author \cite{ManinYano}.

\vspace{.2in}
\noindent{\bf Other Types of Algorithms} We have dealt with classical deterministic algorithms.
Can we do the same things for other types of algorithms. For example, it would be nice to
have universal properties of categories of non-deterministic algorithms, probabilistic algorithms,
parallel algorithms, quantum algorithms, etc. In some sense, with the use of our bracket operation,
we have already dealt with parallel algorithms. Quantum algorithms are a little harder because the
no-cloning theorem does not permit one to have a fully defined product which can lead to
a diagonalization map $x\mapsto (x,x)$.

\vspace{.2in}
\noindent{\bf More Relational Axioms.} It would be interesting to look at other relations that
tell when two programs are essentially the same. With each new relation, we will get different
categories of algorithms and a projection from the old category of algorithms to the new one. With each
new relation, one must find the universal properties of the category of algorithms.

\vspace{.2in}

\noindent{\bf Canonical Presentations of Algorithms.}
Looking at the equivalent trees, one might ask whether there a canonical presentation
of an algorithm. Perhaps we can push up the recursions to the top of the tree, or
perhaps push the brackets to the bottom. This would be most useful for program
correctness and other areas of computer science.

In a sense, Kleene's Theorem on partial recursive functions is an example of
a canonical presentation of an algorithm. It says
that for every computable function, there exists at least one tree-like description of the function
such that the root of the tree is the only minimization in the entire tree.

\vspace{.2in}
\noindent{\bf When are Two Programs Really Different Algorithms.} Is there a way to tell when two programs are
really different algorithms? There is a subbranch of homotopy theory called obstruction theory.
Obstruction theory asks when are two topological spaces in different homotopy classes of spaces.
Is there an obstruction theory of algorithms?

\vspace{.2in}
\noindent{\bf Other Universal Objects in $\catxn$.} We only looked at one element of
$\catxn$ namely $\PRalg$. But there are many other elements that are worthy of study. Given
an arbitrary function $f:\N \raa\N$, consider the category $\C_f$ with $\N$ as its only object and $f$ as its only
non-trivial morphism. The free $\catxn$ category over $\C_f$ is the category of primitive
recursive functions with oracle computations from $f$. It would be nice to frame relative computation theory and
complexity theory from this perspective.

\vspace{.2in}
\noindent{\bf Proof Theory.} There are many similarities between our work and work in proof theory. Many times, one
sees two proofs that are essentially the same. In a sense, Lambek and Scott's excellent book \cite{Lambekandscott}
has the proof theory version of this paper. They look at equivalence classes of proofs to get categories with extra structure.
There is a belief that a program/algorithm implementing a function $f$ is a proof of the fact that $f(x)=y$. Following this
intuition, there should be a very strong relationship between our work and the work done in proof theory. It would be nice to
formalize this relationship. The work of Maietti (e.g. \cite{Maietti}) is in this direction.

\vspace{.2in}
\noindent{\bf A Language Independent Definition of Algorithms.} Our definition of algorithm is dependent
on the language of primitive recursive functions. We could have, no doubt, done the
same thing for other languages. The intuitive notion of an algorithm is language independent.
Can we find a definition of an algorithm that does not depend on
any language?

Consider the set of {\it all} programs in {\it all} programming
languages. Call this set {\bf Programs}. Partition  this set by the
different programming languages that make the programs. So there
will be a subset of {\bf Programs} called {\bf Java}, a subset
called {\bf C++}, and a subset {\bf PL/1} etc. There is also a
subset called {\bf Primitive Recursive} which will contain all the
trees that we discussed in Section 3. There will be functions
between these different subsets. We might call these functions
(non-optimizing) {\it compilers}. They take as input a program from
one programming language and output a program in another programming
language. In some sense {\bf Primitive Recursive} is initial for all
the these sets. By initial we mean that there are compilers going
out of it. There are few compilers going into it. The reason for
this is that in C++ one can program the Ackermann function. One can
not do this in {\bf Primitive Recursive}. (There are, of course,
weaker programming languages than primitive recursive functions, but
we ignore them here.)

For each subset of programs, e.g. {\bf Progs1}, there is a an equivalence
relation $\approx_{\bf Progs1}$ or $\approx_1$ that tells when two
programs in the subset are essentially the same. If $C$ is a compiler from {\bf Progs1} to {\bf Progs2} then we demand
that if two programs in {\bf Progs1} are essentially the same, then the compiled versions
of each of these programs will also be essentially the same, i.e., for any two programs $P$ and $P'$ in {\bf Progs1},
$$P \approx_1 P' \qquad \Rightarrow \qquad C(P) \approx_2 C(P').$$
We also demand that if there are two compilers, then the two compiled programs will be essentially the same,
$$\mbox{ For all programs } P, \qquad C(P) \approx_2 C'(P).$$

Now place the following equivalence relation $\equiv$ on the set {\bf Programs} of {\it all} programs. Two programs
are equivalent if they are the in the same programming language and they are essentially the same, i.e.,
$$P\equiv P'\mbox{ if there exists a relation } \approx_i \mbox{ such that } P \approx_i P'$$
and two programs are equivalent if they are in different programming languages but there exists a compiler that
takes one to the other,
$$ P \equiv P' \mbox{ if there exists a compiler }C \mbox{ and } C(P)=P'.$$

We have now placed an equivalence relation on the set of all programs that tells when two programs are
essentially the same. The equivalence classes of {\bf Programs/$\equiv$} are algorithms. This definition
does not depend on any preferred programming languages. There is much work to do in order to formulate these
ideas correctly. It would also be nice to list the properties of {\bf Algorithms = Programs/$\equiv$}.


\begin{thebibliography}{99}
\bibitem{Abramsky} S. Abramsky. ``Temperley-Lieb algebra: from knot theory to logic and computation via quantum mechanics''.
In {\it Mathematics of Quantum Computing and Technology},
Goong Chen, Louis Kauffman and Sam Lomonaco, eds. Taylor and Francis, 515--558, 2007.
\bibitem{BaezStay}J.C. Baez and M. Stay. ``Physics, Topology, Logic and Computation:
A Rosetta Stone'' http://math.ucr.edu/home/baez/rosetta.pdf. Downloaded from the web on November 30, 2009.
\bibitem{TTT1} M. Barr and C. Wells.
{\it Toposes, triples and theories.} Grundlehren der Mathematischen
Wissenschaften, 278. Springer-Verlag, New York, (1985).
\bibitem{TTT} M. Barr and C. Wells.
{\it Category Theory for Computing Science.} Prentice Hall (1990).
\bibitem{BlassGur}A. Blass, Y. Gurevich. ``Algorithms: A Quest for Absolute Definitions.'' Available
at http://research.microsoft.com/en-us/um/people/gurevich/opera/164.pdf. Downloaded February 5, 2009.
\bibitem{BlassDerGur} A. Blass, N. Dershowitz, and Y. Gurevich. ``When are two algorithms the same?''. Available
at http://arxiv.org/PS\_cache/arxiv/pdf/0811/0811.0811v1.pdf. Downloaded Feburary 5, 2009.
\bibitem{Burroni}A. Burroni. Récursivité graphique. I. Catégorie des fonctions récursives primitives formelles.
[Graph recursiveness. I. The category of formal primitive recursive functions] {\it Cahiers
TopologieGéom. Différentielle Catég}, 27, 49–79, 1986.
\bibitem{Clote} P. Clote. Computational Models and Function Algebras. Handbook of Computability Theory.
Volume 140 (Studies in Logic and the Foundations of Mathematics)
North Holland; 1 edition (October 15, 1999) ISBN-10: 0444898824
\bibitem{Corman} T.H. Corman, C.E. Leiserson, R.L. Rivest, and C. Stein;
{\it Introduction to Algorithms, Second Edition.}
McGraw-Hill (2002).
\bibitem{Dean} W. Dean. {\it What algorithms could not be.} 2006 Thesis in Department of Philosophy. Rutgers University.
\bibitem{Gladstone} M. D. Gladstone. ``A Reduction of the Recursion Scheme''. J. of Symbolic Logic, 32, 4,. 505-508 (1967).
\bibitem{Johnstone}P. T. Johnstone. {\it Sketches of an Elephant: A Topos Theory Compendium: Volume 2}. OUP,
2002.
\bibitem{Knuth} D.E. Knuth. {\it The Art of Computer Programing: Volume 1 / Fundamental Algorithms.} Third Edition.
Addison-Wesley. 1997.
\bibitem{Knuth2} D.E. Knuth. {\it Selected Papers on Computer Science.} Cambridge University Press. 1996.
\bibitem{Lambekandscott} J. Lambek and P.J. Scott. {\it Introduction to higher order categorical logic.} Cambridge University
Press. 1986.
\bibitem{MacLane} Saunders Mac Lane. {\it Categories for the Working Mathematician,} Second Edition. Springer, 1998.
\bibitem{Maietti} M.E. Maietti ``Joyal's Arithmetic Universe via Type Theory''. Electronic notes in Theoretical Computer Science.
69 (2003). Published by Elsevier Science B. V.
\bibitem{Manin} Yu. I. Manin. {\it A course in Mathematical Logic for Mathematicians-- Second Edition}. Springer. October 2009.
\bibitem{ManinYano} Yu. I. Manin, and N.S. Yanofsky. ``Notes on the Recursive Operad''. Work in Progress.
\bibitem{Mosch} Y.N. Moschovakis. ``What Is an Algorithm?'' Available on web. URL: http://citeseerx.ist.psu.edu/viewdoc/summary?doi=10.1.1.7.5576.
(Last accessed on April 3, 2009.)
\bibitem{Pareandroman} R. Pare and L. Roman. ``Monoidal Categories with Natural Numbers Object.'' Studia Logica XLVIII,3 (1989), 361-376.
\bibitem{Roman} L. Roman. ``Cartesian Categories with Natural Numbers Object.'' Journal of Pure and Applied Algebra. 58 (1989), 267-278.
\bibitem{mft} M.-F. Thibault. ``Prerecursive categories'', Journal of Pure and Applied Algebra,. vol. 24,(1982), 79–93.
\bibitem{Yano} N. Yanofsky. ``Towards a Definition of an Algorithm''. Available on the web at
http://arxiv.org/PS\_cache/math/pdf/0602/0602053v2.pdf. (Last accessed on April 3, 2009.)
\end{thebibliography}
\end{document}